\documentclass[a4paper,10pt]{amsart}

\usepackage{amssymb}
\usepackage[OT4]{fontenc}
\usepackage{texdraw}
\usepackage{fullpage}

\theoremstyle{definition}
\newtheorem{definition}{Definition}

\theoremstyle{plain}

\theoremstyle{remark}

\newtheorem{example}[definition]{Example}

\DeclareMathOperator{\edim}{edim}

\DeclareMathOperator{\red}{red}

\DeclareMathOperator{\rev}{rev}

\DeclareMathOperator{\leng}{length}
\DeclareMathOperator{\cut}{cut}
\DeclareMathOperator{\cutr}{\overline{cut}}

\def\ass{\longleftarrow}

\def\sys{\mathcal{L}}

\def\diag{\mathfrak{diag}}

\def\symbx{\mathbf{x}}

\def\aif{{\bf if} }
\def\then{{\bf then} }
\def\return{{\bf return} }
\def\afor{{\bf for} }
\def\aelse{{\bf else} }
\def\endif{{\bf end if} }
\def\endfor{{\bf end for} }
\def\arepeat{{\bf repeat} }
\def\endrepeat{{\bf end repeat} }
\def\ado{{\bf do} }
\def\foreach{{\bf for each} }
\def\endforeach{{\bf end for each} }
\def\achoose{{\bf choose} }
\def\auntil{{\bf until} }

\def\Algreduce{{\sc reduce}}
\def\Algseqreduce{{\sc sequence-reduce}}
\def\Algred{{\sc red}}
\def\Algredout{{\sc redout}}
\def\Algtopreduce{{\sc top-reduce}}
\def\Alghtails{{\sc h-tails}}
\def\Alghtailsall{{\sc ltails}}
\def\Algaddhtailsall{{\sc atails}}
\def\Algsymbred{{\sc symb-reduce}}
\def\Algtails{{\sc tails}}
\def\Algsetbign{{\sc setbign}}
\def\Algsetbigndt{{\sc setbign23}}
\def\Algsetbignb{{\sc setbignb}}
\def\Algsetnb{{\sc setnb}}
\def\Algsetnba{{\sc setnba}}
\def\Algsetpb{{\sc setpb}}
\def\Algsetpba{{\sc setpba}}
\def\Algmatrix{{\sc i-matrix}}
\def\Algrank{{\sc rank}}
\def\Algnonspec{{\sc ns}}
\def\Algnonspecdiags{{\sc check}}
\def\Algcheck{{\sc ch}}
\def\Algsystem{{\sc finalnba}}
\def\Algline{{\sc take-line}}
\def\Algcremona{{\sc cremona}}
\def\Algsort{{\sc sort}}
\def\Algspec{{\sc spec}}

\let\to\longrightarrow

\newenvironment{algorithm}[1]{\medskip\noindent{\bf Algorithm} #1\\}{\medskip}

\begin{document}

\title[Special homogeneous linear systems...]{Special homogeneous linear systems on Hirzebruch surfaces --- algorithmic issues}

\author{Marcin Dumnicki}

\dedicatory{
Institute of Mathematics, Jagiellonian University, \\
ul. \L{}ojasiewicza 6, 30-348 Krak\'ow, Poland \\
Email address: Marcin.Dumnicki@im.uj.edu.pl\\
}



\begin{abstract}
We present algorithms used in the computational part of
the article ``Special homogeneous linear systems on Hirzebruch surfaces''.
\end{abstract}

\maketitle

The main aim of this paper is to give the detailed description
of algorithms used in proving the main Theorem 6 in \cite{mainp}.
We will not repeat the definitions presented in \cite{mainp}.

All algorithms are presented in self-explaining pseudo-code.
We will use indentation to make our algorithms easier to read.
The $\ass$ means an assignment, i.e. in the line\\
$A \ass B$;\\
we force $A$ to be equal to $B$. We use only integers,
so $m \geq 2$ means that $m \in \{2,3,4,\dots\}$.
The command {\bf return} finishes our algorithm immediately,
so further commands (if any) won't be executed.
The control structures ({\bf if}, {\bf for}, {\bf repeat} and so on)
will be used in two versions. The first one, with only one command
executed:

\noindent
\aif $\dots$ \then command;\\
and the second one, with possibility of more that one command to execute:

\noindent
\aif \dots \then\\
\hspace*{0.5cm} command A;\\
\hspace*{0.5cm} command B;\\
\hspace*{0.5cm} $\dots$;\\
\endif

\section{Basic algorithms}

The first algorithm is used to $m$-reduce a given diagram,
according to \cite[Definition 20]{mainp}.

\begin{algorithm}{\Algreduce}

\noindent 
\begin{tabular}{rl}
{\bf Input:} & $m \geq 2$, a diagram $D=\diag(a_1,\dots,a_k).$ \\
{\bf Output:} & $\diag(b_1,\dots,b_k)=\red_m(D)$
or {\sc not reducible} if $D$ is not $m$-reducible.
\end{tabular}
\\

\noindent
\aif $k < m$ \then \return {\sc not reducible};\\
\afor $j=1,\dots,k$ \ado $b_j \ass a_j$;\\
$U \ass \varnothing$;\\
\afor $j=k,\dots,k-m+1$ \ado\\
\hspace*{0.5cm} \aif $a_j < m$ \then $r \ass a_j$ \aelse $r \ass \max (\{1,\dots,m\} \setminus U)$;\\
\hspace*{0.5cm} $b_j \ass a_j - r$;\\
\hspace*{0.5cm} $U \ass U \cup \{r\}$;\\
\endfor\\
\aif $U = \{1,\dots,m\}$ \then \return $\diag(b_1,\dots,b_k)$ \aelse \return {\sc not reducible};\\
\end{algorithm}

\begin{example}
Let us compute $\text{\Algreduce}(3,\diag(5,5,4,2))$.
First we put $(b_1,b_2,b_3,b_4)=(5,5,4,2)$ and take $U = \varnothing$. Next, for
$j=4,3,2$ we will proceed as in the ``{\bf for}'' loop.
For $j=4$ we have $a_4 = 2 < 3 = m$, so we take $r=2$, $b_4=2-2=0$, and
we add $\{2\}$ to the set $U$.
The second step, for $j=3$, gives $r=\max(\{1,2,3\}\setminus \{2\})=3$,
so $r=3$ and $b_3=1$. Before passing to the third step, we put $U=\{2,3\}$.
In the last step we will have
$r=\max(\{1,2,3\}\setminus\{2,3\})=1$, so $b_2=4$.
At the end we have $U=\{1,2,3\}$, so
$\text{\Algreduce}(3,\diag(5,5,4,2))=\diag(5,4,1)$.
\end{example}

We will apply a sequence of reductions to a diagram $D$, so we define
an auxiliary algorithm \Algseqreduce.
By $\red_m^{(k)}$ we denote
$$\red_m^{(k)} = \underbrace{\red_m \circ \dots \circ \red_m}_{k}.$$

\begin{algorithm}{\Algseqreduce}

\noindent 
\begin{tabular}{rl}
{\bf Input:} & $m \geq 2$, a number $k \geq 1$, a diagram $D$.\\
{\bf Output:} & $\red_m^{(k)}(D)$ or {\sc not reducible}
if reduction fails at some step.
\end{tabular}
\\

\noindent
\arepeat $k$ times\\
\hspace*{0.5cm} $D \ass \text{\Algreduce}(m,D)$;\\
\hspace*{0.5cm} \aif $D = \text{{\sc not reducible}}$ \then \return {\sc not reducible};\\
\endrepeat\\
\return $D$;\\
\end{algorithm}

\begin{example}
Let us compute $\text{\Algseqreduce}(4,3,\diag([4]^{\times 3},[5]^{\times 5}))$.
The reductions goes as follows:
$$
\diag([4]^{\times 3},5,5,5,5,5) \to \diag([4]^{\times 3},5,4,3,2,1) \to \diag(4,4,4,5) \to \diag(3,2,1,1),$$
so 
$$\text{\Algseqreduce}(4,3,\diag([4]^{\times 3},[5]^{\times 5})) = \diag(3,2,1,1),$$
while
$$\text{\Algseqreduce}(4,4,\diag([4]^{\times 3},[5]^{\times 5})) = \text{{\sc not reducible}}.$$
\end{example}

The following two algorithms will be used in the algorithm {\Algcheck}.
The first one simply reduces all diagrams from a given set (ommiting not
reducible ones), the second finds all diagrams from a given set $\mathcal D$
which cannot be reduced to some diagram from the second given set $\mathcal G$.

\begin{algorithm}{\Algred}

\noindent 
\begin{tabular}{rl}
{\bf Input:} & $m \geq 2$, a number $k \geq 1$, a set $\mathcal D$ of diagrams.\\
{\bf Output:} & the set $\mathcal G = \{ \red_m^{(k)}(D) : D \in \mathcal D, D \text{ is $m$-reducible $k$ times}\}$.
\end{tabular}
\\

\noindent
$\mathcal G \ass \varnothing$;\\
\foreach $D \in \mathcal D$ \ado\\
\hspace*{0.5cm} $D \ass \text{\Algseqreduce}(m,k,D)$;\\
\hspace*{0.5cm} \aif $D \neq \text{{\sc not reducible}}$ \then $\mathcal G \ass \mathcal G \cup \{D\}$;\\
\endforeach\\
\return $\mathcal G$;
\end{algorithm}

\begin{algorithm}{\Algredout}

\noindent 
\begin{tabular}{rl}
{\bf Input:} & $m \geq 2$, a number $k \geq 1$, a set $\mathcal D$ of diagrams, a set $\mathcal G$ of diagrams.\\
{\bf Output:} & the set $\{ D \in \mathcal D : \red_m^{(k)}(D) \notin \mathcal G\}$.
\end{tabular}
\\

\noindent
\foreach $D \in \mathcal D$ \ado\\
\hspace*{0.5cm} $G \ass \text{\Algseqreduce}(m,k,D)$;\\
\hspace*{0.5cm} \aif $G \in \mathcal G$ \then $\mathcal D \ass \mathcal D \setminus \{D\}$;\\
\endforeach\\
\return $\mathcal D$;
\end{algorithm}

We will often reduce a diagram $D$ as many times as possible, so
we define an auxiliary algorithm \Algtopreduce.

\begin{algorithm}{\Algtopreduce}

\noindent 
\begin{tabular}{rl}
{\bf Input:} & $m \geq 2$, a diagram $D=\diag(a_1,\dots,a_k)$. \\
{\bf Output:} & $G=\diag(b_1,\dots,b_k)=\red_m(\red_m(\dots (D) \dots ))$
such that $G$ is \\
& not $m$-reducible.
\end{tabular}
\\

\noindent
\arepeat\\
\hspace*{0.5cm} $G \ass \text{\Algreduce}(m,D)$;\\
\hspace*{0.5cm} \aif $G=\text{{\sc not reducible}}$ \then \return $D$;\\
\hspace*{0.5cm} $D \ass G$;\\
\endrepeat\\
\end{algorithm}

Observe that if $D$ is not $m$-reducible then $\text{\Algtopreduce}(D)=D$.

\begin{example}
Let us compute $\text{\Algtopreduce}(3,\diag(5,5,4,2))$.
We have 
\begin{align*}
\text{\Algreduce}(3,\diag(5,5,4,2)) & = \diag(5,4,1), \\
\text{\Algreduce}(3,\diag(5,4,1)) & =\diag(3,1) \\
\text{\Algreduce}(3,\diag(3,1)) & =\text{{\sc not reducible}}.
\end{align*}
Hence $\text{\Algtopreduce}(3,\diag(5,5,4,2))=\diag(3,1)$.
\end{example}

Now we present an algorithm to find all admissible
$h$-$\diag(b_1,\dots,b_{m-1})$-tails for multiplicity $m$, see \cite[Definition 35]{mainp}.
By the length of a diagram $D=\diag(a_1,\dots,a_k)$ we denote the number of
its non-zero layers,
$$\leng(D) = \# \{j : a_j > 0\}.$$

\begin{algorithm}{\Alghtails}

\noindent 
\begin{tabular}{rl}
{\bf Input:} & $m \geq 2$, $h > m$, $D=\diag(b_1,\dots,b_{m-1})$. \\
{\bf Output:} & the set $\mathcal D$ of all admissible $h$-$D$-tails
for multiplicity $m$, \\
& or {\sc error} if some reduction stops too early.
\end{tabular}
\\

\noindent
$\mathcal D \ass \varnothing$;\\
\arepeat\\
\hspace*{0.5cm} $\mathcal D \ass \mathcal D \cup \{D\}$;\\
\hspace*{0.5cm} $D \ass \diag(h) + D$;\\
\hspace*{0.5cm} $D \ass \text{\Algtopreduce}(m,D)$;\\
\hspace*{0.5cm} \aif $\leng(D) \geq m$ \then \return {\sc error};\\
\hspace*{0.5cm} \aif $D \in \mathcal D$ \then \return $\mathcal D$;\\
\endrepeat\\
\end{algorithm}

\begin{example}
We will find all admissible $4$-$\diag(0,0)$-tails for multiplicity 3.
The set $\mathcal D$ is empty at the beginning, and we put
$D = \diag(0,0) = \varnothing$ into $\mathcal D$.
Now we take new $D = \diag(4)$, reduce it as many times as possible,
but, in fact, $\text{\Algtopreduce}(3,\diag(4))=\diag(4)$.
So we go at the beginning of the ``{\bf repeat}'' loop and add
$\diag(4)$ to the set $\mathcal D$. Now we take new $D=\diag(4)+\diag(4)=\diag(4,4)$,
still it cannot be $3$-reduced. Thus we add it into
$\mathcal D$, which now is equal to $\{\varnothing, \diag(4), \diag(4,4)\}$.
Taking $D=\diag(4)+\diag(4,4)=\diag(4,4,4)$ we obtain
$\text{\Algtopreduce}(3,\diag(4,4,4))=\varnothing$. Since
the last diagram already belongs to $\mathcal D$, the algorithm
terminates.

Observe that the size of $\mathcal D$ depends also on $h$ and $\diag(b_1,\dots,b_{m-1})$.
For example,
\begin{align*}
\text{\Alghtails}(3,5,\diag(0,0))= \{ &
\varnothing,\diag(5),\diag(5,5),\diag(3),\diag(5,3), \\
& \diag(4,3),\diag(4,2),\diag(4,1),\diag(3,1)\},
\end{align*}
while
$$\# \text{\Alghtails}(3,5,\diag(1,0)) = \# \text{\Alghtails}(3,5,\diag(0,0)) + 4.$$
\end{example}

\begin{example}
\label{exhtails}
We will compute $\text{\Alghtails}(5,9,\diag(0,0,0,0))$.
The computations can be written in the following short form:
\begin{gather*}
\varnothing \to \diag(9) \to \diag(9,9) \to \diag(9,9,9) \to \diag(9,9,9,9) \to \\
\diag(9,9,9,9,9) \to \diag(8,7,6,5,4) \to \diag(7,5,3) \to \diag(9,7,5,3) \to \\
\diag(9,9,7,5,3) \to \diag(8,7,3) \to \diag(9,8,7,3) \to \diag(9,9,8,7,3) \to \\
\diag(8,7,4,2) \to \diag(9,8,7,4,2) \to \diag(8,5,2) \to \diag(9,8,5,2) \to \\
\diag(9,9,8,5,2) \to \diag(8,6,4) \to \diag(9,8,6,4) \to \diag(9,9,8,6,4) \to \\
\diag(8,7,5,1) \to \diag(9,8,7,5,1) \to \diag(7,5,3)
\end{gather*}
and we finish, since the last diagram has been found earlier.
\end{example}

Since we are interested in collecting all admissible $h$-$D$-tails
for all $D \in \mathcal D$, we present an auxiliary algorithm,
called \Alghtailsall.

\begin{algorithm}{\Alghtailsall}

\noindent 
\begin{tabular}{rl}
{\bf Input:} & $m \geq 2$, $h > m$, a set $\mathcal D$ of diagrams.\\
{\bf Output:} & the set $\mathcal G$ of all admissible $h$-$D$-tails
for multiplicity $m$ and $D \in \mathcal D$,\\
&  or {\sc error} if some reduction stops too early.
\end{tabular}
\\

\noindent
$\mathcal G \ass \varnothing$;\\
\foreach $D \in \mathcal D$ \ado\\
\hspace*{0.5cm} $\mathcal H \ass \text{\Alghtails}(m,h,D)$;\\
\hspace*{0.5cm} \aif $\mathcal H=\text{{\sc error}}$ \then \return {\sc error};\\
\hspace*{0.5cm} $\mathcal G \ass \mathcal G \cup \mathcal H$;\\
\endforeach\\
\return $\mathcal G$;\\
\end{algorithm}

Sometimes we also want to find all (top)reductions of diagrams
of the form
$$\diag([h]^{\times n}) + D : D \in \mathcal D$$
for some fixed $n$ and a set $\mathcal D$ of diagrams.

\begin{algorithm}{\Algaddhtailsall}

\noindent 
\begin{tabular}{rl}
{\bf Input:} & $m \geq 2$, $h > m$, $n > 0$, a set $\mathcal D$ of diagrams.\\
{\bf Output:} & the set $\mathcal G$ of all diagram from $\{\diag([h]^{\times n}) + D : D \in \mathcal D\}$ \\
& reduced as many times as possible, \\
& or {\sc error} if some reduction stops too early.
\end{tabular}
\\

\noindent
$\mathcal G \ass \varnothing$;\\
\foreach $D \in \mathcal D$ \ado\\
\hspace*{0.5cm} $G \ass \diag([h]^{\times n}) + D$;\\
\hspace*{0.5cm} $G \ass \text{\Algtopreduce}(m,G)$;\\
\hspace*{0.5cm} \aif $\leng(G) \geq m$ \then \return {\sc error};\\
\hspace*{0.5cm} $\mathcal G \ass \mathcal G \cup \{G\}$;\\
\endforeach\\
\return $\mathcal G$;\\
\end{algorithm}

\begin{example}
Let us compute $\text{\Algaddhtailsall}(3,5,2,\mathcal D)$
for
$$\mathcal D = \{\varnothing, \diag(5),\diag(5,5),\diag(3)\}.$$
For each $D \in \mathcal D$ we execute
$\text{\Algtopreduce}(3,\diag(5,5)+D)$.
The result is
$$\{\diag(5,5),\diag(3),\diag(5,3),\diag(4,3)\}.$$
Observe that
$\text{\Algaddhtailsall}(3,4,2,\mathcal D)=\text{{\sc error}}$,
since 
$$\text{\Algtopreduce}(3,\diag(4,4,5))=\diag(4,2,2),$$
which is too long.
\end{example}

Now our aim is to enumerate all admissible $\diag(a_1,\dots,a_m)$-tails
for multiplicity $m$, see \cite[Definition 39]{mainp}. All admissible
tails could be found by iterating symbolic reductions, see
\cite[Definition 37]{mainp} and the discussion after Example 40. Therefore we
present an algorithm to produce all symbolic reductions of a given
diagram of the form $\diag(a_1,\dots,a_m,[\symbx]^{\times k})$. This amounts
to substitute $[\symbx]^{\times k}$ by all reasonable integers and reduce
obtained diagrams. Observe that for $k=0$, the symbolic reduction
is equal to the $m$-reduction.

If $D = \diag(a_1,\dots,a_k)$ then by $\cut(D, r)$ we denote
the diagram given by
$$\cut(D,r) = \begin{cases}
\diag(a_1,\dots,a_r) & r \leq k, \\
D & r \geq k.
\end{cases}$$

\begin{algorithm}{\Algsymbred}

\noindent 
\begin{tabular}{rl}
{\bf Input:} & $m \geq 2$, $\diag(a_1,\dots,a_m,[\symbx]^{\times k})$, $0 \leq k \leq m-1$.\\
{\bf Output:} & the set $\mathcal D$ of all symbolic reductions of $\diag(a_1,\dots,a_m,[\symbx]^{\times k})$.
\end{tabular}
\\

\noindent
$\mathcal D \ass \varnothing$;\\
\aif $k=0$ \then\\
\hspace*{0.5cm} $D \ass \text{\Algreduce}(m,\diag(a_1,\dots,a_m))$;\\
\hspace*{0.5cm} \aif $D \neq \text{{\sc not reducible}}$ \then $\mathcal D \ass \{D\}$;\\
\hspace*{0.5cm} \return $\mathcal D$;\\
\endif\\
\foreach $(c_1,\dots,c_k)$ satisfying $\min\{m+1,a_m\}\geq c_1 \geq c_2 \geq \dots \geq c_k$ \ado\\
\hspace*{0.5cm} $D \ass \diag(a_1,\dots,a_m,c_1,\dots,c_k)$;\\
\hspace*{0.5cm} $D \ass \text{\Algreduce}(m,D)$;\\
\hspace*{0.5cm} \aif $D \neq \text{{\sc not reducible}}$ \then\\
\hspace*{1.0cm} $G \ass \cut(D,m)$;\\
\hspace*{1.0cm} $\ell \ass \leng(D) - \leng(G)$;\\
\hspace*{1.0cm} $G \ass G + \diag([\symbx]^{\times \ell})$;\\
\hspace*{1.0cm} $\mathcal D \ass \mathcal D \cup \{G\}$;\\
\hspace*{0.5cm} \endif\\
\endforeach\\
\return $\mathcal D$;\\
\end{algorithm}

\begin{example}
\label{exsymbred}
Let us compute $\text{\Algsymbred}(3,\diag(5,5,5,\symbx,\symbx))$.
In the ``{\bf foreach}'' loop we must consider all pairs $(c_1,c_2)$ satisfying
$4 \geq c_1 \geq c_2$. For each such pair we take
$\diag(5,5,5,c_1,c_2)$, $3$-reduce it, and change the fourth and fifth number into
symbols $\symbx$, if necessary. We present the computations in the
following table.
$$
\begin{array}{c|c|c|c}
(c_1,c_2) & \diag(5,5,5,c_1,c_2) & \red_3(\diag(5,5,5,c_1,c_2)) & \text{the result} \\
\hline
(4,4) & \diag(5,5,5,4,4) & \diag(5,5,4,2,1) & \diag(5,5,4,\symbx,\symbx) \\
(4,3) & \diag(5,5,5,4,3) & \diag(5,5,4,2) & \diag(5,5,4,\symbx) \\
(4,2) & \diag(5,5,5,4,2) & \diag(5,5,4,1) & \diag(5,5,4,\symbx) \\
(4,1) & \diag(5,5,5,4,1) & \diag(5,5,3,1) & \diag(5,5,3,\symbx) \\
(4,0) & \diag(5,5,5,4) & \diag(5,4,3,1) & \diag(5,4,3,\symbx) \\
(3,3) & \diag(5,5,5,3,3) & \diag(5,5,4,1) & \diag(5,5,4,\symbx) \\
(3,2) & \diag(5,5,5,3,2) & \diag(5,5,4) & \diag(5,5,4) \\
(3,1) & \diag(5,5,5,3,1) & \diag(5,5,3) & \diag(5,5,3) \\
(3,0) & \diag(5,5,5,3) & \diag(5,4,3) & \diag(5,4,3) \\
(2,2) & \diag(5,5,5,2,2) & \text{{\sc not reducible}} & \\
(2,1) & \diag(5,5,5,2,1) & \diag(5,5,2) & \diag(5,5,2) \\
(2,0) & \diag(5,5,5,2) & \diag(5,4,2) & \diag(5,4,2) \\
(1,1) & \diag(5,5,5,1,1) & \text{{\sc not reducible}} & \\
(1,0) & \diag(5,5,5,1) & \diag(5,3,2) & \diag(5,3,2) \\
(0,0) & \diag(5,5,5) & \diag(4,3,2) & \diag(4,3,2)
\end{array}
$$
\end{example}

Now we can enumerate all admissible $\diag(a_1,\dots,a_m)$-tails
for multiplicity $m$. To do this, we will consider all symbolic reductions
of symbolic reductions of $\dots$ and so on. Each diagram
obtained in this way, which is short enough (i.e. without $\symbx$ and with
length at most $m-1$), satisfies the desired property.

\begin{algorithm}{\Algtails}

\noindent 
\begin{tabular}{rl}
{\bf Input:} & $m \geq 2$, $\diag(a_1,\dots,a_m)$.\\
{\bf Output:} & the set $\mathcal D$ of all admissible $\diag(a_1,\dots,a_m)$-tails.
\end{tabular}
\\

\noindent
$\mathcal D \ass \varnothing$;\\
$\mathcal W \ass \{\diag(a_1,\dots,a_m,[\symbx]^{\times (m-1)})\}$;\\
\arepeat\\
\hspace*{0.5cm} \achoose $W \in \mathcal W$;\\
\hspace*{0.5cm} $\mathcal W \ass \mathcal W \setminus \{W\}$;\\
\hspace*{0.5cm} $\mathcal R \ass \text{\Algsymbred}(m,W)$;\\
\hspace*{0.5cm} $\mathcal W \ass \mathcal W \cup \mathcal R$;\\
\hspace*{0.5cm} \foreach $D \in \mathcal R$ \ado\\
\hspace*{1.0cm} \aif $\leng(D) < m$ \then $\mathcal D \ass \mathcal D \cup \{D\}$;\\
\hspace*{0.5cm} \endforeach\\
\auntil $\mathcal W = \varnothing$\\
\end{algorithm}

To show that the above algorithm terminates after a finite number of steps,
observe that in each step, after choosing $W \in \mathcal W$ and
producing $\mathcal R = \text{\Algsymbred}(m,W)$, we have
$$\# W > \max \{ \# R : R \in \mathcal R\},$$
where
$$\# \diag(a_1,\dots,a_k,[\symbx]^{\ell}) = a_1 + \dots + a_k.$$

\begin{example}
\label{ex255}
We will find $\text{\Algtails}(2,\diag(5,5))$. The example for $m=3$ would
be a bit too long.
The idea is to compute consecutive symbolic reductions
of $\diag(5,5,\symbx)$.
In the first step we proceed as in Example \ref{exsymbred}.
$$
\begin{array}{c|c|c|c|c}
D & c_1 & D(c_1) & \red_2(D(c_1)) & \text{the result} \\ \hline
\diag(5,5,\symbx) & 3 & \diag(5,5,3) & \diag(5,4,1) & \diag(5,4,\symbx) \\
\diag(5,5,\symbx) & 2 & \diag(5,5,2) & \diag(5,4) & \diag(5,4) \\
\diag(5,5,\symbx) & 1 & \diag(5,5,1) & \diag(5,3) & \diag(5,3) \\
\diag(5,5,\symbx) & 0 & \diag(5,5) & \diag(4,3) & \diag(4,3)
\end{array}
$$
In the next step we take all obtained diagrams and perform all possible
symbolic reductions.
$$
\begin{array}{c|c|c|c|c}
D & c_1 & D(c_1) & \red_2(D(c_1)) & \text{the result} \\ \hline
\diag(5,4,\symbx) & 3 & \diag(5,4,3) & \diag(5,3,1) & \diag(5,3,\symbx) \\
\diag(5,4,\symbx) & 2 & \diag(5,4,2) & \diag(5,3) & \diag(5,3) \\
\diag(5,4,\symbx) & 1 & \diag(5,4,1) & \diag(5,2) & \diag(5,2) \\
\diag(5,4,\symbx) & 0 & \diag(5,4) & \diag(4,2) & \diag(4,2) \\
\diag(5,4) & & \diag(5,4) & \diag(4,2) & \diag(4,2) \\
\diag(5,3) & & \diag(5,3) & \diag(4,1) & \diag(4,1) \\
\diag(4,3) & & \diag(4,3) & \diag(3,1) & \diag(3,1)
\end{array}
$$
Again, in the third step:
$$
\begin{array}{c|c|c|c|c}
D & c_1 & D(c_1) & \red_2(D(c_1)) & \text{the result} \\ \hline
\diag(5,3,\symbx) & 3 & \diag(5,3,3) & \diag(5,2,1) & \diag(5,2,\symbx) \\
\diag(5,3,\symbx) & 2 & \diag(5,3,2) & \diag(5,2) & \diag(5,2) \\
\diag(5,3,\symbx) & 1 & \diag(5,3,1) & \diag(5,1) & \diag(5,1) \\
\diag(5,3,\symbx) & 0 & \diag(5,3) & \diag(4,1) & \diag(4,1) \\
\diag(5,3) & & \diag(5,3) & \diag(4,1) & \diag(4,1) \\
\diag(5,2) & & \diag(5,2) & \diag(4) & \diag(4) \\
\diag(4,2) & & \diag(4,2) & \diag(3) & \diag(3) \\
\diag(4,1) & & \diag(4,1) & \diag(2) & \diag(2) \\
\diag(3,1) & & \diag(3,1) & \diag(1) & \diag(1)
\end{array}
$$
The diagrams with length 1 are no more $2$-reducible, so they won't
produce any additional admissible tail.
We present the fourth step:
$$
\begin{array}{c|c|c|c|c}
D & c_1 & D(c_1) & \red_2(D(c_1)) & \text{the result} \\ \hline
\diag(5,2,\symbx) & 2 & \diag(5,2,2) & \diag(5,1) & \diag(5,1) \\
\diag(5,2,\symbx) & 1 & \diag(5,2,1) & \diag(5) & \diag(5) \\
\diag(5,2,\symbx) & 0 & \diag(5,2) & \diag(4) & \diag(4) \\
\diag(5,2) & & \diag(5,2) & \diag(4) & \diag(4) \\
\diag(5,1) & & \diag(5,1) & \diag(3) & \diag(3) \\
\diag(4,1) & & \diag(4,1) & \diag(2) & \diag(2)
\end{array}
$$
In the final step we must reduce $\diag(5,1)$ to obtain $\diag(3)$.
We collect all obtained diagrams of length at most 1, thus
$$
\text{\Algtails}(2,\diag(5,5)) = \{ \diag(5), \diag(4), \diag(3), \diag(2), \diag(1)\}.
$$
Observe that we can skip some of the above reducing.
\end{example}

\section{Algorithms to compute sets $\mathcal D$}

In \cite[Section 7]{mainp} we construct various sets of diagrams.
Each set serves for showing that some given family of systems
contains only non-special ones. To be more precise, we define a family
$\mathcal S$ of systems together with a finite set $\mathcal D$ of diagrams
such that
if for each $D \in \mathcal D$ and
$r = \lfloor \frac{\# D}{\binom{m+1}{2}} \rfloor$ the systems
$\sys(D;m^{\times r})$ and $\sys(D;m^{\times (r+1)})$ are non-special
then $\mathcal S$ contains only non-special systems.

The first algorithm computes the set $\mathcal D$ from \cite[Proposition 42]{mainp}.
For $D=\diag(a_1,\dots,a_k)$ let $\rev(D) = \diag(a_k,\dots,a_1)$.
Similarly, for a set $\mathcal D$ of diagrams, let
$$\rev(\mathcal D) = \{ \rev(D) : D \in \mathcal D\}.$$

\begin{algorithm}{\Algsetbign}

\noindent 
\begin{tabular}{rl}
{\bf Input:} & $m \geq 4$, $N \geq m$.\\
{\bf Output:} & the set $\mathcal D$ from Proposition 42.
\end{tabular}
\\

\noindent
$\mathcal D \ass \varnothing$;\\
$\mathcal L \ass \text{\Alghtails}(m,m+1,\diag([0]^{\times (m-1)}))$;\\
\afor $j=m+2,\dots,2m-3$ \ado\\
\hspace*{0.5cm} $\mathcal L \ass \text{\Algaddhtailsall}(m,j,N,\mathcal L)$;\\
\hspace*{0.5cm} $\mathcal L \ass \text{\Alghtailsall}(m,j,\mathcal L)$;\\
\endfor\\
$\mathcal L \ass \text{\Alghtailsall}(m,2m-2,\mathcal L)$;\\
$\mathcal R \ass \text{\Algtails}(m,\diag([2m-1]^{\times m}))$;\\
\foreach $(L,R) \in \mathcal L \times \mathcal R$ \ado\\
\hspace*{0.5cm} $\mathcal D \ass \mathcal D \cup \{\rev(L) + \diag([2m-2]^{\times N}) + R\}$;\\
\endforeach\\
\return $\mathcal D$;\\
\end{algorithm}

\begin{example}
\label{exsetbign}
We will show the example for $m=5$, $N=11$ (which is a part of our computation
to prove Theorem 6 in \cite{mainp}). We will not present all the details,
since the output would be too big.
In our case we do the following:

\noindent
$\mathcal L \ass \text{\Alghtails}(5,6,\diag(0,0,0,0))$;\\
$\mathcal L \ass \text{\Algaddhtailsall}(5,7,11,\mathcal L)$;\\
$\mathcal L \ass \text{\Alghtailsall}(5,7,\mathcal L)$;\\
$\mathcal L \ass \text{\Alghtailsall}(5,8,\mathcal L)$;\\
$\mathcal R \ass \text{\Algtails}(5,\diag(9,9,9,9,9))$;\\

In the first step we obtain
$$\mathcal L = \{\varnothing,\diag(6),\diag(6,6),\diag(6,6,6),\diag(6,6,6,6)\}.$$
In the second step, for each $D \in \mathcal L$, we take
$\text{\Algtopreduce}(5,\diag([7]^{\times 11})+D)$.
After reducing, we will have
$$\mathcal L = \{\diag(7,6,4),\diag(7,7,6,3),\diag(6,4,3,1),\diag(5),\diag(7,4)\}.$$
In the third step we look for \Algtopreduce{} of all diagrams of the form
$$\diag([7]^{\times k}) + D, \qquad k \geq 0, \: D \in \mathcal L.$$
We will not enumerate all of them, since after this step, $\# \mathcal L = 53$.
In the fourth step we look for
$$\diag([8]^{\times k}) + D, \qquad k \geq 0, \: D \in \mathcal L,$$
and the resulting set contains $119$ diagrams.
Now we look for admissible tails. After computations, we obtain
$\mathcal R$ of cardinality $147$.

Now we must ``glue'' diagrams from $\mathcal L$ and $\mathcal R$ to produce
$119 \cdot 147 = 17493$ diagrams in $\mathcal D$.
One can check that, for example,
$$\diag(8,6,3,1) \in \mathcal L, \quad \diag(7,6,5,4) \in \mathcal R,$$
so we have
$$\diag(1,3,6,8,[8]^{\times 11},7,6,5,4) \in \mathcal D.$$
\end{example}

Our next algorithm computes the set $\mathcal D$ from Proposition 44 in \cite{mainp}.

\begin{algorithm}{\Algsetbigndt}

\noindent 
\begin{tabular}{rl}
{\bf Input:} & $2 \leq m \leq 3$, $N \geq m$.\\
{\bf Output:} & the set $\mathcal D$ from Proposition 44.
\end{tabular}
\\

\noindent
$\mathcal D \ass \varnothing$;\\
$\mathcal L \ass \text{\Alghtails}(m,m+1,\diag([0]^{\times (m-1)}))$;\\
$\mathcal L \ass \text{\Alghtailsall}(m,m+2,\mathcal L)$;\\
$\mathcal R \ass \text{\Algtails}(m,\diag([m+3]^{\times m}))$;\\
\foreach $(L,R) \in \mathcal L \times \mathcal R$ \ado\\
\hspace*{0.5cm} $\mathcal D \ass \mathcal D \cup \{\rev(L) + \diag([m+2]^{\times N}) + R\}$;\\
\endforeach\\
\return $\mathcal D$;\\
\end{algorithm}

\begin{example}
The example for $m=3$ would be very nice and illustrating,
but also too long. So we will
compute $\text{\Algsetbigndt}(2,2)$.
We have three steps:

\noindent
$\mathcal L \ass \text{\Alghtails}(2,3,\diag(0))$;\\
$\mathcal L \ass \text{\Alghtailsall}(2,4,\mathcal L)$;\\
$\mathcal R \ass \text{\Algtails}(2,\diag(5,5))$;\\

In the first step we obtain
$$\mathcal L = \{\varnothing, \diag(3)\}.$$
In the next step
$$\mathcal L = \{\diag(4),\diag(3),\diag(2),\diag(1)\}.$$
In Example \ref{ex255} we have shown that
$$\mathcal R = \{\diag(5),\diag(4),\diag(3),\diag(2),\diag(1)\}.$$
So the final set is
$$\mathcal D = \{\diag(a,[4]^{\times 2},b) : 1\leq a \leq 4, 1\leq b \leq 5\}.$$
\end{example}

Our next algorithm computes the set $\mathcal D$ from
\cite[Proposition 47]{mainp}.

\begin{algorithm}{\Algsetbignb}

\noindent 
\begin{tabular}{rl}
{\bf Input:} & $m \geq 2$, $N \geq m$, $b \geq m+2$.\\
{\bf Output:} & the set $\mathcal D$ from Proposition 47.
\end{tabular}
\\

\noindent
$\mathcal D \ass \varnothing$;\\
$\mathcal L \ass \text{\Alghtails}(m,m+1,\diag([0]^{\times (m-1)}))$;\\
\afor $j=m+2,\dots,b-1$ \ado\\
\hspace*{0.5cm} $\mathcal L \ass \text{\Algaddhtailsall}(m,j,N,\mathcal L)$;\\
\hspace*{0.5cm} $\mathcal L \ass \text{\Alghtailsall}(m,j,\mathcal L)$;\\
\endfor\\
$\mathcal L \ass \text{\Alghtailsall}(m,b,\mathcal L)$;\\
$\mathcal R \ass \text{\Alghtails}(m,b+1,\diag([0]^{\times (m-1)}))$;\\
\foreach $(L,R) \in \mathcal L \times \mathcal R$ \ado\\
\hspace*{0.5cm} $\mathcal D \ass \mathcal D \cup \{\rev(L) + \diag([b]^{\times N}) + R\}$;\\
\endforeach\\
\return $\mathcal D$;\\
\end{algorithm}

\begin{example}
We will show the example for $m=5$, $N=11$, $b=8$ (which is a part of our computation
to prove Theorem 6 in \cite{mainp}).
In our case we do the following:

\noindent
$\mathcal L \ass \text{\Alghtails}(5,6,\diag(0,0,0,0))$;\\
$\mathcal L \ass \text{\Algaddhtailsall}(5,7,11,\mathcal L)$;\\
$\mathcal L \ass \text{\Alghtailsall}(5,7,\mathcal L)$;\\
$\mathcal L \ass \text{\Alghtailsall}(5,8,\mathcal L)$;\\
$\mathcal R \ass \text{\Alghtails}(5,9,\diag(0,0,0,0))$;\\

This is very similar to what we did in Example \ref{exsetbign}, except
for the last step.
So we will have $\# \mathcal L = 119$. Now we must look for all possible
top-reductions of $\diag([9]^{\times k})$. By Example \ref{exhtails}
we have $\# \mathcal R = 15$, so
$\# \mathcal D = 15 \cdot 119 = 1785$.
\end{example}

The next algorithm computes the set $\mathcal D$ from
\cite[Proposition 50]{mainp}. For a diagram
$D=\diag(a_1,\dots,a_k)$ let
$$\cutr(D,\ell) = \begin{cases}
\diag(a_{k-\ell+1},a_{k-\ell+2},\dots,a_k) & \ell \leq k, \\
D & \ell > k.
\end{cases}$$

\begin{algorithm}{\Algsetnb}

\noindent 
\begin{tabular}{rl}
{\bf Input:} & $m \geq 2$, $n \geq 2$, $B \geq 2m-1$.\\
{\bf Output:} & the set $\mathcal D$ from Proposition 50.
\end{tabular}
\\

\noindent
$\mathcal D \ass \varnothing$;\\
$G \ass \diag([m+1]^{\times n},\dots,[B]^{\times n},B+1)$;\\
$H \ass \cutr(G,m)$;\\
$K \ass \cut(G,n(B-m)-m+1)$;\\
$\mathcal R \ass \text{\Algtails}(m,H)$;\\
\foreach $R \in \mathcal R$ \ado\\
\hspace*{0.5cm} $\mathcal D \ass \mathcal D \cup \{K + R\}$;\\
\endforeach\\
\return $\mathcal D$;\\
\end{algorithm}

\begin{example}
Let us compute $\text{\Algsetnb}(3,2,6)$.
We will have
$$G = \diag(4,4,5,5,6,6,7),$$
so we take
$$H = \diag(6,6,7), \qquad K = \diag(4,4,5,5).$$
Now we must find $\mathcal R = \text{\Algtails}(m,H)$ by considering
all symbolic reductions of $\diag(6,6,7,\symbx,\symbx)$, and take
$\mathcal D = \{ \diag(4,4,5,5) + R : R \in \mathcal R \}$.
\end{example}

The next algorithm computes the set $\mathcal D$ from
\cite[Proposition 52]{mainp}.

\begin{algorithm}{\Algsetnba}

\noindent 
\begin{tabular}{rl}
{\bf Input:} & $m \geq 2$, $n \geq 2$, $b \geq m+1$, $A \geq 0$.\\
{\bf Output:} & the set $\mathcal D$ from Proposition 52.
\end{tabular}
\\

\noindent
$\mathcal D \ass \varnothing$;\\
$\mathcal R \ass \text{\Alghtails}(m,b+1,\diag([0]^{\times (m-1)}))$;\\
\foreach $R \in \mathcal R$ \ado\\
\hspace*{0.5cm} $\mathcal D \ass \mathcal D \cup \{\diag([m+1]^{\times n},\dots,[b]^{\times n},[b+1]^{\times (A+1)}) + R\}$;\\
\endforeach\\
\return $\mathcal D$;\\
\end{algorithm}

\begin{example}
For $(m,n,b,A)=(3,2,5,0)$ we will have
\begin{align*}
\mathcal R & = \text{\Alghtails}(3,6,\diag(0,0)) \\
& = \{ \varnothing, \diag(6), \diag(6,6), \diag(4,2), \diag(5,1) \},
\end{align*}
so
$$D = \{ \diag(4,4,5,5,6) + R : R \in \mathcal R\}.$$
\end{example}

The next algorithm computes the set $\mathcal D$ from
\cite[Proposition 54]{mainp}.

\begin{algorithm}{\Algsetpb}

\noindent 
\begin{tabular}{rl}
{\bf Input:} & $m \geq 2$, $B \geq 3(m-1)$.\\
{\bf Output:} & the set $\mathcal D$ from Proposition 54.
\end{tabular}
\\

\noindent
$\mathcal D \ass \varnothing$;\\
$\mathcal R \ass \text{\Algtails}(m,\diag(B-m+2,B-m+3,\dots,B+1))$;\\
\foreach $R \in \mathcal R$ \ado\\
\hspace*{0.5cm} $\mathcal D \ass \mathcal D \cup \{\diag(1,2,\dots,B-m+1) + R\}$;\\
\endforeach\\
\return $\mathcal D$;\\
\end{algorithm}

\begin{example}
We will compute $\text{\Algsetpb}(3,9)$.
Hence our computation starts with
$$\mathcal R = \text{\Algtails}(3,\diag(8,9,10)).$$
We obtain $\# \mathcal R = 28$ and take
$28$ diagrams of the form
$$\diag(1,2,\dots,7) + R : R \in \mathcal R.$$
For example, we will have $\diag(1,2,3,4,5,6,7,6,5) \in \mathcal D$.
\end{example}

The next algorithm computes the set $\mathcal D$ from
\cite[Proposition 56]{mainp}.

\begin{algorithm}{\Algsetpba}

\noindent 
\begin{tabular}{rl}
{\bf Input:} & $m \geq 2$, $B \geq m$, $A \geq b$.\\
{\bf Output:} & the set $\mathcal D$ from Proposition 56.
\end{tabular}
\\

\noindent
$\mathcal D \ass \varnothing$;\\
$\mathcal R \ass \text{\Alghtails}(m,b+1,\diag([0]^{\times (m-1)}))$;\\
\foreach $R \in \mathcal R$ \ado\\
\hspace*{0.5cm} $\mathcal D \ass \mathcal D \cup \{\diag([b+1]^{\times (A+1)}) + R\}$;\\
\endforeach\\
\return $\mathcal D$;\\
\end{algorithm}

\begin{example}
\label{exsetpba}
It is easy to check that
\begin{align*}
\text{\Algsetpba}(3,7,7) = \{ & \: \diag([8]^{\times 8}), \diag([8]^{\times 9}), \diag([8]^{\times 10}), \diag([8]^{\times 8},5,1), \\
& \: \diag([8]^{\times 8},6,2), \diag([8]^{\times 8},7,3), \diag([8]^{\times 8},7,5)\}.
\end{align*}
\end{example}

\section{Checking non-speciality}

We begin with an algorithm to decide whether a given system
$\sys(D;m^{\times r})$ is special or not. The computations
will be performed over $\mathbb F_p$
and for some randomly chosen coordinates of points. Therefore, if
our ``specialized'' system is non-special then obviously the general one
is also non-special. In the opposite case we only know that our
method does not work. We begin with preparing the matrix for our system
(see \cite{jsc}).

\begin{algorithm}{\Algmatrix}

\noindent 
\begin{tabular}{rl}
{\bf Input:} & a diagram $D$, $m\geq 2$, $r \geq 1$, $p_1,\dots,p_r \in \mathbb F_p^2$.\\
{\bf Output:} & the matrix $M$ associated to $\sys(D;mp_1,\dots,mp_r)$.\\
{\bf Property:} & computations over $\mathbb F_p$.
\end{tabular}
\\

\noindent
$\mathcal M \ass \{x^{\alpha}y^{\beta} : (\alpha,\beta) \in D\}$;\\
$f_{\mathcal M} \ass \text{ a one-to-one correspondence from } \mathcal M \text{ to } \{1,\dots,\# D\}$;\\
$\mathcal C \ass \{ (k,d_x,d_y) : 1 \leq k \leq r, d_x+d_y<m \}$;\\
$f_{\mathcal C} \ass \text{ a one-to-one correspondence from } \mathcal C \text{ to } \{1,\dots,r\binom{m+1}{2}\}$;\\
\foreach $(x^\alpha y^\beta,(k,d_x,d_y)) \in \mathcal M \times \mathcal C$ \ado\\
\hspace*{0.5cm} $g \ass \frac{\partial^{(d_x+d_y)}}{\partial x^{d_x} \partial y^{d_y}}(x^\alpha y^\beta)$;\\
\hspace*{0.5cm} $g \ass g(p_k)$;\\
\hspace*{0.5cm} $M(f_{\mathcal M}(x^\alpha y^\beta),f_{\mathcal C}(k,d_x,d_y)) \ass g$;\\
\endforeach\\
\return $M$;\\
\end{algorithm}

\begin{example}
Let us compute $M$ for $D=\diag(3,2,1)$, $m^{\times r}=2^{\times 2}$, $p_1=(0,0)$, $p_2=(2,1)$.
The set of monomials 
$$\mathcal M = \{ 1,x,y,x^2,xy,y^2 \}$$
will be ordered by $f_{\mathcal M}$ as above, the set of
conditions will also be ordered as indicated:
$$\mathcal C = \{(1,0,0),(1,1,0),(1,0,1),(2,0,0),(2,1,0),(2,0,1)\}.$$
Now, for example, take $x^2$ and $(1,1,0)$. The polynomial
$$g = \frac{\partial \: x^2}{\partial x}(0,0) = 0$$
will be inserted into $M$ in the $f_{\mathcal M}(x^2)=4$th row and
the $f_{\mathcal C}(1,1,0)=2$nd column.
For the same monomial and condition $(2,1,0)$ we will have
$$g = \frac{\partial \: x^2}{\partial x}(2,1) = 4$$
inserted into $M[4,5]$.
The entire matrix is equal to
$$
M = \left[
\begin{array}{cccccc}
1 & 0 & 0 & 1 & 0 & 0 \\
0 & 1 & 0 & 2 & 1 & 0 \\
0 & 0 & 1 & 1 & 0 & 1 \\
0 & 0 & 0 & 4 & 4 & 0 \\
0 & 0 & 0 & 2 & 1 & 2 \\
0 & 0 & 0 & 1 & 0 & 1
\end{array}
\right].
$$
\end{example}

For a matrix $M$ over $\mathbb F_p$, let
$\text{\Algrank}(M)$ denote the rank of $M$ computed, for example,
by using the Gauss elimination.

\begin{algorithm}{\Algnonspec}

\noindent 
\begin{tabular}{rl}
{\bf Input:} & $m\geq 2$, $r \geq 1$, a diagram $D$, a number of tries $t \geq 1$.\\
{\bf Output:} & {\sc non-special} or {\sc not decided}, \\
              & {\sc non-special} implies that $\sys(D;m^{\times r})$ is non-special.\\
{\bf Property:} & computations over $\mathbb F_p$.
\end{tabular}
\\

\noindent
\arepeat $t$ times\\
\hspace*{0.5cm} $(p_1,\dots,p_r) \ass $ randomly chosen points in $\mathbb F_p^2$;\\
\hspace*{0.5cm} $M \ass \text{\Algmatrix}(D,m,r,p_1,\dots,p_r)$;\\
\hspace*{0.5cm} $k \ass \text{\Algrank}(M)$;\\
\hspace*{0.5cm} \aif $k=\min\{\#D,r\binom{m+1}{2}\}$ \then \return {\sc non-special};\\
\endrepeat\\
\return {\sc not decided};\\
\end{algorithm}

\begin{example}
We will have
$\text{\Algnonspec}(2,2,\diag(3,2,1),t) = \text{{\sc not decided}}$,
since the system $\sys(\diag(3,2,1);2,2)=\sys(2;2,2)$ is special.

Taking $D=\diag(2,1)$, $m=1$ and $r=3$ we will obtain
$\text{\Algnonspec}(1,3,\diag(2,1),t) = \text{{\sc non-special}}$
if and only if the algorithm chooses $p_1,p_2,p_3$ not lying on a line.
For a non-special system, the result {\sc non-special} is much more probable if the number $t$ of tries is big.
During computations, it appeared that ``big'' in our case means $t \geq 6$.
\end{example}

We will check non-speciality of all systems $\sys(D;m^{\times r})$
for $D \in \mathcal D$, fixed $m$ and all $r \geq 1$.

\begin{algorithm}{\Algnonspecdiags}

\noindent 
\begin{tabular}{rl}
{\bf Input:} & $m\geq 2$, a set of diagrams $\mathcal D$, a number of tries $t \geq 1$.\\
{\bf Output:} & a set $\mathcal G \subset \mathcal D$ such that for every $G \in \mathcal G$, $r \geq 1$, \\
              & $\sys(G;m^{\times r})$ is non-special.
\end{tabular}
\\

\noindent
$\mathcal G \ass \varnothing$;\\
\foreach $D \in \mathcal D$ \ado\\
\hspace*{0.5cm} $r \ass \left\lfloor \frac{\# D}{\binom{m+1}{2}} \right\rfloor$;\\
\hspace*{0.5cm} \aif $\text{\Algnonspec}(m,r,D,t) = \text{{\sc non-special}}$ \then\\
\hspace*{1.0cm} \aif $\text{\Algnonspec}(m,r+1,D,t) = \text{{\sc non-special}}$ \then\\
\hspace*{1.5cm} $\mathcal G \ass \mathcal G \cup \{D\}$;\\
\hspace*{1.0cm} \endif\\
\hspace*{0.5cm} \endif\\
\endforeach\\
\return $\mathcal G$;\\
\end{algorithm}

Observe that the above algorithm is sufficient to check whether all
diagrams in $\mathcal D$ gives non-special systems for a fixed multiplicity.
However, running it on the set $\mathcal D$ (of cardinality $17493$)
from Example \ref{exsetbign} would consume too much time.
Therefore we will reduce all diagrams from $\mathcal D$ several times
(this should decrease the number of diagrams) and check whether they are non-special.
I yes, we are done due to \cite[Theorem 27]{mainp}. If no,
we must deal with diagrams that reduces to special ones.
This will be explained in more details after presenting the algorithm.

\begin{algorithm}{\Algcheck}

\noindent 
\begin{tabular}{rl}
{\bf Input:} & $m\geq 2$, a set of diagrams $\mathcal D$,\\
             & a number $u \geq 0$ of reductions, \\
             & a number $v \geq 0$ of reductions performed on reversed diagrams. \\
{\bf Output:} & {\sc ok} or {\sc not decided}, \\
              & {\sc ok} implies that for every $D \in \mathcal D$, $r \geq 1$, $\sys(D;m^{\times r})$ is non-special.
\end{tabular}
\\

\noindent
\aif $u > 0$ \then\\
\hspace*{0.5cm} $\mathcal R \ass \text{\Algred}(m,u,\mathcal D)$;\\
\hspace*{0.5cm} $\mathcal R \ass \text{\Algnonspecdiags}(m,\mathcal R,6)$;\\
\hspace*{0.5cm} $\mathcal D \ass \text{\Algredout}(m,u,\mathcal D,\mathcal R)$;\\
\endif\\
\aif $v > 0$ \then\\
\hspace*{0.5cm} $\mathcal D \ass \rev(\mathcal D)$;\\
\hspace*{0.5cm} $\mathcal R \ass \text{\Algred}(m,v,\mathcal D)$;\\
\hspace*{0.5cm} $\mathcal R \ass \text{\Algnonspecdiags}(m,\mathcal R,6)$;\\
\hspace*{0.5cm} $\mathcal D \ass \text{\Algredout}(m,v,\mathcal D,\mathcal R)$;\\
\endif\\
$\mathcal R \ass \text{\Algnonspecdiags}(m,\mathcal D,16)$;\\
\aif $\mathcal R = \mathcal D$ \then \return {\sc ok} \aelse \return {\sc not decided};\\
\end{algorithm}

\begin{example}
Let us deal with $\text{\Algcheck}(3,\mathcal D,8,0)$
for $\mathcal D$ from Example \ref{exsetpba}.
\begin{align*}
\mathcal D = \{ & \: \diag([8]^{\times 8}), \diag([8]^{\times 9}), \diag([8]^{\times 10}), \diag([8]^{\times 8},5,1), \\
& \: \diag([8]^{\times 8},6,2), \diag([8]^{\times 8},7,3), \diag([8]^{\times 8},7,5)\}.
\end{align*}
Each diagram must be $3$-reduced 8 times.
$$
\begin{array}{c|c}
D & \red_{3}^{(8)}(D) \\
\hline
\diag([8]^{\times 8}) & \diag(8,6,2) \\
\diag([8]^{\times 9}) & \diag(8,8,6,2) \\
\diag([8]^{\times 10}) & \diag(8,8,8,6,2) \\
\diag([8]^{\times 8},5,1) & \diag(8,7,5,2) \\
\diag([8]^{\times 8},6,2) & \diag(8,8,6,2) \\
\diag([8]^{\times 8},7,3) & \diag(8,8,7,3) \\
\diag([8]^{\times 8},7,5) & \diag(8,8,7,5)
\end{array}
$$
We end up with the set $\mathcal R$ containing 6 diagrams (reducing decreased the
number of cases). We can check that
$\mathcal R = \text{\Algnonspecdiags}(3,\mathcal R,6)$,
so we are done and the result is {\sc ok}. An additional advantage lies
in the size of matrices, since each $m$-reduction decreases the number of
rows and columns by $\binom{m+1}{2}$.
\end{example}

\begin{example}
The set $\mathcal D$ from Example \ref{exsetbign} contains $17493$ diagrams.
We will run $\text{\Algcheck}(5,\mathcal D,3,0)$.
So we must $5$-reduce every diagram in $\mathcal D$ three times,
which gives the set $\mathcal R$ of cardinality $6234$.
All of these diagrams appeared to be non-special, so
$$\mathcal R = \text{\Algnonspecdiags}(5,\mathcal R,6)$$
and we are done.
\end{example}

\begin{example}
We will present the number of diagrams involved in computing
$\text{\Algcheck}(6,\mathcal D,14,13)$ for
$$\mathcal D = \text{\Algsetbignb}(6,51,8).$$
The set $\mathcal D$ contains $5472$ diagrams.
$4617$ of them can be $6$-reduced $14$ times and we obtain
the set $\mathcal R$ of $2991$ diagrams.
By checking speciality we obtain that $2832$ diagrams from $\mathcal R$
are non-special, while the rest is probably special.
So in $\mathcal D$ we have $855$ not-reducible diagrams together with
$250$ that reduces to special ones.
Now we reverse $1105$ diagrams, reduce them (all are reducible) $13$ times to
obtain the set with $562$ diagrams. Again not all of them are non-special,
we are left with $46$ diagrams that reduces to $24$ special ones.
\end{example}

In the next algorithm we deal with systems $\sys_n(a,b;m^{\times r})$
for fixed $m$, $n$, $a$ and $b$. Our aim is to identify
those $r$, for which the system is special.

\begin{algorithm}{\Algsystem}

\noindent 
\begin{tabular}{rl}
{\bf Input:} & $m\geq 2$, $n \geq 0$, $a,b \geq 0$.\\
{\bf Output:} & the set $\mathcal L \subset \mathbb N$ such that \\
              & if $r \notin \mathcal L$ then $\sys_n(a,b;m^{\times r})$ is non-special.
\end{tabular}
\\

\noindent
$\mathcal L \ass \varnothing$;\\
\aif $n=0$ \then $D \ass \diag([a+1]^{\times (b+1)})$;\\
\aif $n\geq 2$ \then $D \ass \diag([1]^{\times n},\dots,[b]^{\times n},[b+1]^{\times (a+1)})$;\\
$r \ass \left\lfloor \frac{\# D}{\binom{m+1}{2}} \right\rfloor$;\\
\arepeat\\
\hspace*{0.5cm} $A \ass \text{\Algnonspec}(m,r,D,16)$;\\
\hspace*{0.5cm} \aif $A = \text{{\sc not decided}}$ \then $\mathcal L \ass \mathcal L \cup \{r\}$;\\
\hspace*{0.5cm} $r \ass r-1$;\\
\auntil $A=\text{{\sc non-special}}$;\\
$r \ass \left\lfloor \frac{\# D}{\binom{m+1}{2}} \right\rfloor + 1$;\\
\arepeat\\
\hspace*{0.5cm} $A \ass \text{\Algnonspec}(m,r,D,16)$;\\
\hspace*{0.5cm} \aif $A = \text{{\sc not decided}}$ \then $\mathcal L \ass \mathcal L \cup \{r\}$;\\
\hspace*{0.5cm} $r \ass r+1$;\\
\auntil $A=\text{{\sc non-special}}$;\\
\return $\mathcal L$;
\end{algorithm}

\begin{example}
Let us compute $\text{\Algsystem}(3,0,5,4)$.
We have $D = \diag([6]^{\times 5})$, so $\# D = 30$ and, at the beginning,
$r = 5$.
In the first step we have
$$\text{\Algnonspec}(3,5,D,16) = \text{{\sc not decided}},$$
since in fact $\sys_0(5,4;3^{\times 5})$ is special.
So we take $\mathcal L = \{5\}$ and compute
$$\text{\Algnonspec}(3,6,D,16) = \text{{\sc non-special}}.$$
Since also
$$\text{\Algnonspec}(3,4,D,16) = \text{{\sc non-special}},$$
we finish with $\mathcal L = \{5\}$.
\end{example}

The last group of algorithms checks whether a given system
is $-1$-special, see \cite[Definition 3]{mainp}.
We begin with auxiliary algorithms {\Algline} and {\Algcremona}
(see \cite[Theorem 3]{jsc}).
We put
\begin{align*}
\text{\Algline}(\sys(d;m_1,\dots,m_r)) & = \sys(d-1;m_1-1,m_2-1,m_3,\dots,m_r), \\
\text{\Algcremona}(\sys(d;m_1,\dots,m_r)) & = \sys(d+k;m_1+k,m_2+k,m_3+k,m_4,\dots,m_r)
\end{align*}
for $k=d-m_1-m_2-m_3$.
We will also use $\text{\Algsort}(\sys(d;m_1,\dots,m_r))$ to sort multiplicities
in non-increasing order.

\begin{algorithm}{\Algspec}

\noindent 
\begin{tabular}{rl}
{\bf Input:} & $m\geq 2$, $n \geq 0$, $a,b \geq 0$, $r \geq 1$.\\
{\bf Output:} & {\sc -1-special} if $\sys_n(a,b;m^{\times r})$ is $-1$-special, \\
              & {\sc error} otherwise. \\
{\bf Remark:} & for a system of curves $L$ we define $L_d$ to be the degree,\\
              & $L_{m_j}$ to be the $j$-th multiplicity.
\end{tabular}
\\

\noindent
$d \ass (n+1)b+a$;\\
$m_0 \ass nb+a$;\\
$e \ass \edim \sys_n(a,b;m^{\times r})$;\\
\afor $t=0,\dots,b$ \ado\\
\hspace*{0.5cm} $L \ass \sys(d-t;m_0,m^{\times r},(b-t)^{\times (n+1)})$;\\
\hspace*{0.5cm} $L \ass \text{\Algsort}(L)$;\\
\hspace*{0.5cm} \arepeat\\
\hspace*{1.0cm} \aif $L_d - L_{m_1} - L_{m_2} < 0$ \then\\
\hspace*{1.5cm} $L \ass \text{\Algline}(L)$;\\
\hspace*{1.0cm} \aelse\\
\hspace*{1.5cm} \aif $L_d - L_{m_1} - L_{m_2} - L_{m_3} < 0$ \then\\
\hspace*{2.0cm} $L \ass \text{\Algcremona}(L)$;\\
\hspace*{1.5cm} \endif\\
\hspace*{1.0cm} \endif\\
\hspace*{1.0cm} $L \ass \text{\Algsort}(L)$;\\
\hspace*{0.5cm} \auntil $L_d-L_{m_1}-L_{m_2}-L_{m_3} \geq 0$;\\
\hspace*{0.5cm} \aif $\edim L > e$ \then \return {\sc -1-special};\\
\endfor\\
\return {\sc error};\\
\end{algorithm}

\begin{example}
Let us compute $\text{\Algspec}(6,8,2,8,15)$, so we must consider
$\sys_8(2,8;6^{\times 15})$.
We begin with planar system $L=\sys(74;66,6^{\times 15},\overline{8^{\times 9}})$
with $\edim L = -1$.
For $t=0$ take the system
$\sys(74;66,6^{\times 15},8^{\times 9})$, sort multiplicities to
obtain
$\sys(74;66,8^{\times 9},6^{\times 15})$. Then use {\Algcremona}
four times to obtain $\sys(42;34,8,6^{\times 15})$. Again, use {\Algcremona}
to obtain $\sys(36;28,6^{\times 15},2)$. By the sequence of {\Algcremona}
we transform our system into $\sys(8;2^{\times 15})$. Since
$\edim \sys(8;2^{\times 15}) = -1$, we pass to the case $t=1$.
Now we begin with $\sys(73;66,6^{\times 15},7^{\times 9})$ and, by {\Algcremona},
transform to $\sys(9;6,6,2,1^{\times 13})$. Now we use {\Algline} several times
to produce $\sys(6;3,3,2,1^{\times 13})$. By {\Algcremona} we finish with
$\sys(4;1^{\times 15})$ of negative expected dimension.
For $t=2$ we begin with $\sys(72;66,6^{\times 24})$ and transform it to
$\sys(6;6,6)$. Then, by {\Algline}, we obtain $\sys(0;0)$ of non-negative
expected dimension, so the answer is {\sc -1-special}.
\end{example}

\section{Implementation and results}

All the presented algorithms have been implemented in Free Pascal
and can be downloaded from \cite{www}. They are divided into two kinds,
depending on method of working. The first kind simply works on given
data, the second one prepares batch files with the list of instructions.
For example, the implementation of {\Algred} (``\verb"red.pas"''),
of the first kind, performes sequence reductions on given diagrams (loaded
from the specified file). The algorithm
{\Algsetpb} (of the second kind) prepares the batch file with the following
instructions (for $\text{\Algsetpb}(3,9)$)
{\small
\begin{verbatim}
tails 3 8,9,10,x,x rt
basediag 1 1 7 0 bt
gluediags inempty bt rt diag
\end{verbatim}}
The above instructions run \verb"tails", which produces all
admissible $\diag(8,9,10)$-tails and stores them in the file \verb"rt";
\verb"basediag", which prepares $\diag(1,2,3,4,5,6,7)$ and stores it
in \verb"bt"; \verb"gluediags", which glues diagrams from
\verb"inempty" (by default, it contains only the empty diagram),
\verb"bt" and \verb"rt".

All algorithms with names beginning with {\sc set} are of the second kind, together with
{\Algcheck} and {\Algsystem}. The others are of the first kind.

All algorithms produces \verb"log" files, where the necessary information
is stored. For example, the part of \verb"log" file for multiplicity $2$
contains:
{\small
\begin{verbatim}

XXXXXXXXXXXXXXXXXXXXXXXXXXXXXXXXXXXXXXXXXXXXXXXXXXXXXXXXXX

setbign23 2 2 
result: all systems L_n(a,b)(2^r) are non-special 
for n>=2, a>=0, b>=5, r>=0 

XXXXXXXXXXXXXXXXXXXXXXXXXXXXXXXXXXXXXXXXXXXXXXXXXXXXXXXXXX

ltails (all-h-D-admissible tails) 2 3
tails loaded:

1 tails loaded.
tails found:

3
4 entries used, 2 tails found.
 job finished: 03:11:41:14
*************************************

ltails (all-h-D-admissible tails) 2 4
tails loaded:

3
2 tails loaded.
tails found:

3
4
1
2
11 entries used, 5 tails found.
 job finished: 03:11:41:15
*************************************

basediag 4 0 2 0
base diagram:
4,4
 job finished: 03:11:41:15
*************************************

tails (admissible tails) 2
diagram:
5,5,x
tails found:
4
3
2
1
5
17 entries used, 5 tails found.
 job finished: 03:11:41:15
*************************************

gluediags (glue diagrams)
25 diagrams produced.
 job finished: 03:11:41:15
*************************************

check
multiplicity: 2
diag(4,4,4)  det <> 0 det <> 0
diag(4,4,3)  det <> 0 det <> 0
diag(4,4,2)  det <> 0 det <> 0
diag(4,4,1)  det <> 0 det <> 0
diag(4,4,5)  det <> 0 det <> 0
diag(3,4,4,4)  det <> 0 det <> 0
diag(3,4,4,3)  det <> 0 det <> 0
diag(3,4,4,2)  det <> 0 det <> 0
diag(3,4,4,1)  det <> 0 det <> 0
diag(3,4,4,5)  det <> 0 det <> 0
diag(4,4,4,4)  det <> 0 det <> 0
diag(4,4,4,3)  det <> 0 det <> 0
diag(4,4,4,2)  det <> 0 det <> 0
diag(4,4,4,1)  det <> 0 det <> 0
diag(4,4,4,5)  det <> 0 det <> 0
diag(1,4,4,4)  det <> 0 det <> 0
diag(1,4,4,3)  det <> 0 det <> 0
diag(1,4,4,2)  det <> 0 det <> 0
diag(1,4,4,1)  det <> 0 det <> 0
diag(1,4,4,5)  det <> 0 det <> 0
diag(2,4,4,4)  det <> 0 det <> 0
diag(2,4,4,3)  det <> 0 det <> 0
diag(2,4,4,2)  det <> 0 det <> 0
diag(2,4,4,1)  det <> 0 det <> 0
diag(2,4,4,5)  det <> 0 det <> 0
result: positive.
non-special: 25, special: 0
 job finished: 03:11:41:16
*************************************
\end{verbatim}}

The information which set $\mathcal D$ is considered, is stored in the
preamble. It is then followed by names of programs together with
additional detailed information.
The \verb"shortlog" files contain only preambles and names of
programs, while \verb"infolog" stores only preambles.
The \verb"finitlog" files contain informations on running {\Algspec}.
Each program informs when it has terminated (day:hour:minute:second).

\end{document}